\documentclass[
final
]{dmtcs-episciences}

\usepackage{subfigure}
\usepackage{mathrsfs}

\usepackage{amsmath,amssymb,amsfonts,euscript}
\usepackage{algorithmic, algorithm}
\usepackage{graphicx}
\usepackage[ansinew]{inputenc}
\usepackage{graphicx}
\usepackage{color}
\usepackage{mathrsfs}

\newtheorem{prethm}{{\bf Theorem}}

\newenvironment{thm}{\begin{prethm}{\hspace{-0.5
				em}{\bf}}}{\end{prethm}}

\newtheorem{prepro}{{\bf Theorem}}

\newtheorem{preprop}{{\bf Proposition}}

\newtheorem{precor}{{\bf Corollary}}

\newenvironment{cor}{\begin{precor}{\hspace{-0.5
				em}{\bf}}}{\end{precor}}

\newtheorem{preconj}{{\bf Conjecture}}

\newenvironment{conj}{\begin{preconj}{\hspace{-0.5
				em}{\bf}}}{\end{preconj}}

\newtheorem{predefi}{{\bf Definition}}

\newtheorem{preremark}{{\bf Remark}}

\newenvironment{remark}{\begin{preremark}\rm{\hspace{-0.5
				em}{\bf}}}{\end{preremark}}

\newtheorem{preexample}{{\bf Fact}}

\newenvironment{example}{\begin{preexample}\rm{\hspace{-0.5
				em}{\bf}}}{\end{preexample}}

\newtheorem{prelem}{{\bf Lemma}}

\newtheorem{prelam}{{\bf Lemma}}

\newenvironment{lam}{\begin{prelam}{\hspace{-0.5
				em}{\bf}}}{\end{prelam}}

\newtheorem{preprob}{{\bf Problem}}

\newtheorem{preproof}{{\bf Proof}}

\newtheorem{preali}{{\bf Proof of Theorem 1.}}

\newenvironment{ali}[1]{\begin{preali}{\rm
			#1}\hfill{$\Box$}}{\end{preali}}

\newtheorem{prealii}{{\bf Proof of Theorem 3.}}

\newenvironment{alii}[1]{\begin{prealii}{\rm
			#1}\hfill{$\Box$}}{\end{prealii}}

\newtheorem{prealiii}{{\bf Proof of Theorem 4.}}

\newenvironment{aliii}[1]{\begin{prealiii}{\rm
			#1}\hfill{$\Box$}}{\end{prealiii}}

\newtheorem{prealiiii}{{\bf Proof of Theorem 5.}}

\newenvironment{aliiii}[1]{\begin{prealiiii}{\rm
			#1}\hfill{$\Box$}}{\end{prealiiii}}

\newtheorem{prealiiiii}{{\bf Proof of Theorem 6.}}

\newenvironment{aliiiii}[1]{\begin{prealiiiii}{\rm
			#1}\hfill{$\Box$}}{\end{prealiiiii}}

\author{Ali Dehghan\affiliationmark{1}
  \and  Mohammad-Reza Sadeghi\affiliationmark{1}
  \and  Arash Ahadi\affiliationmark{2}
  \thanks{E-mail Addresses: 
$\mathsf{alidehghan@sce.carleton.ca}$ (Ali Dehghan) $\mathsf{msadeghi@aut.ac.ir}$ (Mohammad-Reza Sadeghi) arash$_{-}$ahadi@mehr.sharif.edu (Arash Ahadi).}}
\title[Sigma Partitioning: Complexity and Random Graphs]{Sigma Partitioning: Complexity and Random Graphs}

\affiliation{
Department of Mathematics and Computer Science, Amirkabir University of Technology, Tehran, Iran\\
Department of Mathematical Sciences, Sharif University of Technology, Tehran, Iran
}
\keywords{Sigma partitioning; Lucky labeling;  Additive coloring; Sigma chromatic number; Computational Complexity; Planar Not-All-Equal 3-SAT; Planar Not-All-Equal 3-SAT Type 2.}
\received{2016-6-18}

\revised{2018-11-13}

\accepted{2018-11-22}
\begin{document}
\publicationdetails{20}{2018}{2}{19}{1534}
\maketitle
\begin{abstract}{
 A $\textit{sigma partitioning}$ of a graph $G$ is a partition of the vertices into sets $P_1, \ldots, P_k$ such that for every two adjacent vertices $u$ and $v$ there is an index $i$ such that $u$ and $v$ have different numbers of neighbors in $P_i$. The $\textit{ sigma number}$ of a graph $G$, denoted by $\sigma(G)$, is the minimum number $k$ such that
	$ G $ has a sigma partitioning $P_1, \ldots, P_k$.
	Also, a $\textit{ lucky labeling}$ of a graph $G$ is a function $ \ell :V(G) \rightarrow \mathbb{N}$, such that for every two adjacent vertices $ v $ and $
	u$ of $ G $, $ \sum_{w \sim v}\ell(w)\neq \sum_{w \sim u}\ell(w) $ ($ x \sim y $
	means that $ x $  and $y$ are adjacent).
	The $\textit{ lucky number}$ of $ G $, denoted by $\eta(G)$, is the minimum number $k $ such that $ G $ has a lucky labeling $ \ell :V(G) \rightarrow \mathbb{N}_k$.
	It was conjectured in [Inform. Process. Lett., 112(4):109--112, 2012] that it is $ \mathbf{NP} $-complete to decide whether $ \eta(G)=2$ for
	a given 3-regular graph $G$. In this work, we prove this conjecture.
	Among other results, 
	we give an upper bound of five for the sigma number of a uniformly random graph.
}\end{abstract}

\section{Introduction}
\label{}
Throughout the paper we denote $\{1,2,\ldots, k\}$ by $\mathbb{N}_k$.
In 2004, Karo\'nski et al. introduced a new coloring of
a graph which is generated via edge labeling \cite{MR2047539}. Let $f : E(G) \rightarrow \mathbb{N}$ be a labeling of the edges of a
graph $G$ by positive integers such that for every two adjacent vertices $v$ and $u$, $\text{Sum}(v) \neq \text{Sum}(u)$, where $\text{Sum}(v)$ denotes the sum of labels of all edges
incident with $v$. It was conjectured in  \cite{MR2047539} that three integer labels
$\mathbb{N}_3$ are sufficient for every connected graph, except
$K_2$. Currently the best bound is five \cite{MR2595676}.
Regarding the computational complexity of this concept, Dudek and Wajc   \cite{David} proved that determining whether a given graph has a
labeling of the edges from $\mathbb{N}_2$ that induces a proper vertex coloring is
$ \mathbf{NP} $-complete. Recently, it was shown that for a given 3-regular graph  $G$ deciding whether $G$ has a
labeling for the edges from $ \{ a,b\}$, ($a\neq b$) that induces a proper vertex coloring is $ \mathbf{NP} $-complete \cite{MR3072733}.

Lucky labeling and sigma partitioning are two vertex
versions of this problem, which were introduced recently by Czerwi\'nski et al. \cite{MR2552893} and Chartrand et al. \cite{MR2729020}.
The
{\it lucky labeling} of a graph $G$ is a function $ \ell :V(G) \rightarrow\mathbb{N}$, such that for every two adjacent vertices $ v $ and $
u$ of $ G $, $ \sum_{w \sim v}\ell(w)\neq \sum_{w \sim u}\ell(w) $ ($ x \sim y $
means that $ x $ and $y$ are adjacent).
The {\it lucky number} of $ G $, denoted by $\eta(G)$, is the minimum number $k $ such that $ G $ has a lucky labeling $ \ell :V(G) \rightarrow \mathbb{N}_k $.
Also, the {\it sigma partitioning} of a graph $G$ is a partition of the vertices into sets $P_1, \ldots, P_k$ such that for every two adjacent vertices $u$ and $v$, there is an index $i$ such that $u$ and $v$ have different numbers of neighbors in $P_i$. The {\it sigma number} of a graph $G$, denoted by $\sigma(G)$, is the minimum number $k$ such that $ G $ has a sigma partitioning $P_1, \ldots, P_k$.

There is an alternative definition for the sigma partitioning which is similar to the definition of lucky labeling.
For a graph $G$, let $ c :V(G) \rightarrow \mathbb{N}$ be a vertex labeling of
$G$. If $k$ labels are used by $c$, then
$c$ is a $k$-labeling of $G$. If for every two adjacent vertices $ v $ and $
u$ of $ G $, $ \sum_{w \sim v}c(w)\neq \sum_{w \sim u}c(w) $, then $c$ is called a  sigma partitioning  of $G$. The minimum number of labels
required in a sigma partitioning is called the   sigma  number  of $G$.
Now, we show that the alternative definition is  equivalent to the first definition. In any sigma partitioning of a graph $G$ with $\sigma(G)$ labels, we can use the set of labels $\{ s^i : 0 \leq i \leq \sigma(G)-1 \}$, where $s$ is a sufficiently large number (it is enough to put $s\geq \Delta(G)+1$).
So the sigma  number is the minimum number $k$ such that the vertices of graph can be partitioned into $k$ sets $P_1, \ldots, P_k$ such that for every edge $uv$, there is an index $i$ that $u$ and $v$ have different numbers of neighbors in $P_i$.

Sigma partitioning and lucky labeling have been studied extensively by several
authors, for instance see \cite{ ahadi, MR3704829, aaa,  MR2729020, MR2552893, MR3589724,  MR2729390, survey}.
Note that lucky labeling is also called an \emph{additive labeling} of a graph \cite{iapprox, additive},
and sigma partitioning is also called the \emph{sigma chromatic number} of a graph \cite{MR2729020}.

The sigma number can also be thought of as a vertex version of the
\emph{detection number} of a graph. For a connected graph $G$ of order $|V(G)| \geq 3$ and a $k$-labeling $c : E(G) \rightarrow \mathbb{N}_k$ of the edges
of $G$, the code of a vertex $v$ of $G$ is the ordered $k$-tuple $(l_1, l_2, \ldots , l_k)$, where $l_i$ is the number of edges incident with $v$ that are labeled $i$. The $k$-labeling $c$ is detectable, if every two adjacent vertices of $G$ have distinct codes. The minimum positive integer $k$ for which $G$ has a detectable $k$-labeling is the {\it detection number} of $G$. It was shown that it is $ \mathbf{NP} $-complete to decide if the detection number of a 3-regular graph is $2$ \cite{havet2012detection}. Also, it was proved that
the detection number of every bipartite graph of minimum degree at least 3 is at most 2 \cite{havet2012detection}.

\begin{remark}
	For a partition $\mathcal{P}= \cup_{i=1}^{k} P_i$ of the vertices of a graph $G$,
	consider the function $f_{\mathcal{P}}:V(G) \rightarrow
	(\mathbb{N} \cup \{0\})^k$ with the map $f_{\mathcal{P}}(x)= (N(x) \cap P_1, \ldots, N(x) \cap P_k)$.
	The function $\mathcal{P}$ is called  sigma partitioning
	if and only if $f_\mathcal{P}$ is a proper coloring for $G$.
	For a given graph $G$, let $\mathcal{P}_k$ be the set of partitions of the vertices such that each partition has at most $k$ parts. Now, consider the following parameter:
	
	$ \Omega_p(G, k)= \max_{\mathcal{P} \in \mathcal{P}_k} \min_ {vu \in E(G)} |f_\mathcal{P}(v)- f_\mathcal{P}(u)|_p,$
	\\ \\
	where $p\geq 1$ and $|.|_p$ is the norm $p$.
	Note that  $p=1$ is the natural case and  $\sigma (G)= \min_k \Omega_p (G, k)\neq 0$ (for all $p$).
\end{remark}

\begin{remark}
	The difference between the sigma number and the  lucky number of a graph can be arbitrarily large. In fact, there are graphs with $\sigma \leq 2$ and
	arbitrary large lucky numbers. For instance, for every $k$ consider a
	complete graph with $k+2 \choose 2$$-1$ vertices
	$\{v_{\alpha\beta}: \alpha \in \mathbb{N}_k \cup \{0\}, \beta \in \mathbb{N}_{\alpha+2}\}$;
	next join each vertex $v_{\alpha\beta}$ to $\alpha$ new isolated
	vertices $v^\gamma_{\alpha\beta}$ for $ \gamma \in \mathbb{N}_{\alpha}$.
	We call this graph  $G$. Note that $G$  has $k(k+1)(2k+7)/6$ leaves.
	We show $\sigma(G)\leq 2$.
	To do this, put
	$A=\{ v_{\alpha\beta} : \beta =1 \} \cup\{v^ \gamma _{\alpha \beta} : \gamma \leq \beta -2\}$
	and let $\overline{A}=V(G) \setminus A$.
	In order to show
	that $(A,\overline{A})$ is a sigma partition one could mention that we only need to consider pairs of
	vertices of the same degree, such as $v_{ia}$ and $v_{ib}$ for some $i$ and $a < b$, in which case $v_{ib}$
	has more neighbors in $A$ than $v_{ia}$ does.
	Thus $(A , \overline{A}) $ is a sigma partitioning.
	On the other hand, let $\ell$
	be a lucky labeling of $G$, for every two adjacent vertices
	$u$, $v\in K_{{k+2 \choose 2}-1}$, we have:

	$$ \displaystyle\sum_{w \sim v \atop w  \text{ is a leaf } } \ell(w) + \ell(u)
	\neq \displaystyle\sum_{w \sim u \atop w \text{ is a leaf } } \ell(w) + \ell(v),$$
	\\ \\
	Therefore, we have:
	
	$$\displaystyle\sum_{w \sim v \atop w \text{ is a leaf } } \ell(w) - \ell(v)
	\neq \displaystyle\sum_{w \sim u \atop w \text{ is a leaf } } \ell(w) - \ell(u).$$
	\\\\
	$\sum_{w \sim v \atop w \text{ is a leaf } } \ell(w) - \ell(v)$
	is one of the $(k+1)\eta(G)-1$ different numbers.
	Thus $$(k+1)\eta(G)- 1 \geq {{k+2} \choose {2}} -2.$$
	Therefore, $\eta(G)\geq k/2$.
\end{remark}

\subsection{Complexity for 3-regular graphs}

For a given graph $G$, we have $\sigma(G)=\eta(G)=1$ if and only if every two adjacent vertices of $G$ have different degrees. We know that $\sigma(G) \leq \chi(G)$ \cite{MR2729020},  on the other hand we have the following conjecture about the lucky number.

\begin{conj}
	{\bf[Additive Coloring Conjecture \rm \cite{MR2552893}]\bf} For every
	graph $G$, $\eta(G)\leq \chi(G)$.
\end{conj}

It is not known 
whether this conjecture is true even for bipartite graphs.
Moreover, it is not even known if $\eta(G)$ is
bounded for bipartite graphs. Recently, Grytczuk et al.  \cite{additive} proved that $\eta(G) \leq 468$ for every
planar graph $G$.

It was shown in \cite{ahadi} that it is $ \mathbf{NP} $-complete to decide for a given planar 3-colorable graph $G$, whether $ \eta(G)=2$. Here, we are interested in the following conjecture:

\begin{conj}\label{C2}
	{\bf \rm \cite{ahadi}\bf}
	It is $ \mathbf{NP} $-complete to decide whether $ \eta(G)=2$ for
	a given 3-regular graph $G$.
\end{conj}

In this work, we prove Conjecture \ref{C2}. In order to prove Conjecture \ref{C2}; first, we prove Theorem \ref{T1}. After that by this theorem we prove the conjecture.

\begin{thm}\label{T1}
	It is $ \mathbf{NP} $-complete to decide  for a given 3-regular graph $G$, whether $ \sigma(G)=2$.
\end{thm}

It is easy to check that for every regular graph $G$, $\sigma(G)=2$ if and only if $\eta(G)=2$. Therefore from Theorem \ref{T1} we have Theorem \ref{NT2} and we can prove Conjecture 2.

\begin{thm}\label{NT2}
	It is $ \mathbf{NP} $-complete to decide  for a given 3-regular graph $G$, whether $ \eta(G)=2$.
\end{thm}

In the proof of Theorem \ref{T1}, for a given formula $ \Psi $, we transform $ \Psi $ into a 3-regular graph $G_{\Psi}$
such that $ \sigma(G_{\Psi})=2$ if and only if $\Psi $ has a Not-All-Equal (NAE) truth assignment, where a NAE assignment is an assignment such that each clause has at least one true literal and at least one false literal. On the other hand, $G_{\Psi}$ has a triangle and is a 3-regular graph, so we have $\chi(G_{\Psi})=3$. Also, since $G_{\Psi}$ is a 3-regular graph, $\sigma(G_{\Psi})\geq 2$.
Thus, $\sigma(G)=\chi(G)$ if and only if $\Psi $ does not have any NAE truth assignment.
Consequently, by Theorem \ref{T1}, we have the following corollary.

\begin{cor}
	It is $\mathbf{NP}$-complete to determine whether $\sigma(G)=\chi(G)$, for a given 3-regular graph $G$.
\end{cor}

It was proven in  \cite{ahadi} that  for every $k\geq 2$, it is
{\bf NP}-complete to decide whether $\eta(G)=k$ for a given graph
$G$. Here, we present the following result for sigma partitioning.

\begin{thm}
	For every $k\geq 3$, it is
	{\bf NP}-complete to decide whether $\sigma(G)= k$ for a given graph
	$G$.
\end{thm}

Consider the problem of
partitioning the vertices of $G$ into $\sigma(G)$ parts, such that this partitioning is a sigma partitioning and some parts have the smallest possible size. We show that for a given planar $3$-regular graph with sigma number two
the following problem is $ \mathbf{NP} $-complete.
\\ \\
{\em  Problem $ \Xi$.}\\
\textsc{Instance}: A planar $3$-regular graph $G$ with  $\sigma(G) = 2$ and a real number $0 < r  <  1$.\\
\textsc{Question}: Does $G$ have a sigma partitioning $c$ with $\sigma(G)$ parts, such that there is a part with at most $r|V(G)|$ vertices?\\

\begin{thm}\label{T4}
	Problem $\Xi$ is $ \mathbf{NP} $-complete.
\end{thm}

\subsection{Some upper bounds}

It was shown that $\sigma(G)= \mathcal{O} (\Delta^2)$ for every graph $G$ \cite{survey}. We prove the following upper bound for triangle-free graphs.

\begin{thm}\label{T6}
	Let $\mathcal{G} = \{G_i\}_{i \in \mathbb{N}}$ be a sequence  of triangle-free graphs. Then for this family,\\ (i) if $\delta =  \Omega(\Delta)$ then $\sigma = \mathcal{O} (1)$ for all graphs in $\mathcal{G}$, except a finite number of them;\\ (ii) if $\delta \geq  3.501 (\ln \Delta)$, then $\sigma = \mathcal{O} (\delta) $.
\end{thm}

Next, we show that almost all graphs have a small sigma number.
Let $G(n,p)$ be
the notation of  Erd\"{o}s-R\'{e}nyi
random graph on $n$ vertices.

\begin{thm}
	For every constant $p$, $0<p<1$, $\sigma(G(n,p))\leq 5$.
\end{thm}

\subsection{Notation}

We conclude the section by fixing some notation which is not defined here. When we say that $f$ if a sigma partitioning for a graph $G$, we mean  that $f:V(G) \rightarrow \mathbb{N}$ such that for every two adjacent vertices $ v $ and $
u$ of $ G $, $ \sum_{w \sim v}f(w)\neq \sum_{w \sim u}f(w) $.
For a vertex $v$ of $G$, let $N(v)$ denote the neighborhood of
$v$ (the set of vertices adjacent to $v$). Let $ N[v]= N(v)\cup \{ v \}$ denote the closed
neighborhood of $v$. Also, for every $v\in V (G)$, $d_G(v)$ denotes the degree of $v$ (for   simplicity we denote $d_G(v)$ by $d(v)$). For a natural number $k$, a graph $G$ is
called a {\it $k$-regular graph } if $d(v) = k$, for each $v \in V(G)$. We denote the maximum degree and the minimum degree
of $G$ by $\Delta(G)$ and $\delta(G)$, respectively.
For $ k\in \mathbb{N} $, a {\it proper vertex $k$-coloring} of $G$ is a function $c:
V(G)\rightarrow \mathbb{N}_k$, such that if $u,v\in V(G)$ are adjacent,
then $c(u)$ and $c(v)$ are different.
The smallest integer $k$ such that
$G$ has a proper vertex $k$-coloring is called the {\it chromatic number} of $G$ and denoted by $\chi(G)$.
We
follow \cite{MR1367739} for terminology and notation which are not defined here.

\section{Complexity results}
\begin{ali}{
		Here, we prove that  for a given 3-regular graph $G$, it is $ \mathbf{NP} $-complete to decide  whether $ \sigma(G)=2$.
		Let $\Psi$ be a $3$-SAT formula with clauses $C=\lbrace
		c_1, \ldots ,c_k\rbrace $ and variables
		$X=\lbrace x_1, \ldots ,x_n\rbrace $. The following problem is $ \mathbf{NP}$-complete \cite{MR1567289}.
		\\ \\
		{\em Not-All-Equal (NAE) 3-SAT .}\\
		\textsc{Instance}: Set $X$ of variables, collection $C$ of clauses over $X$ such that each
		clause $c \in C$ has $|c | = 3$.\\
		\textsc{Question}: Is there a truth assignment for $X$ such that each clause in $C$ has at
		least one true literal and at least one false literal?\\

		\begin{figure}[ht]
			\begin{center}
				\includegraphics[scale=.6]{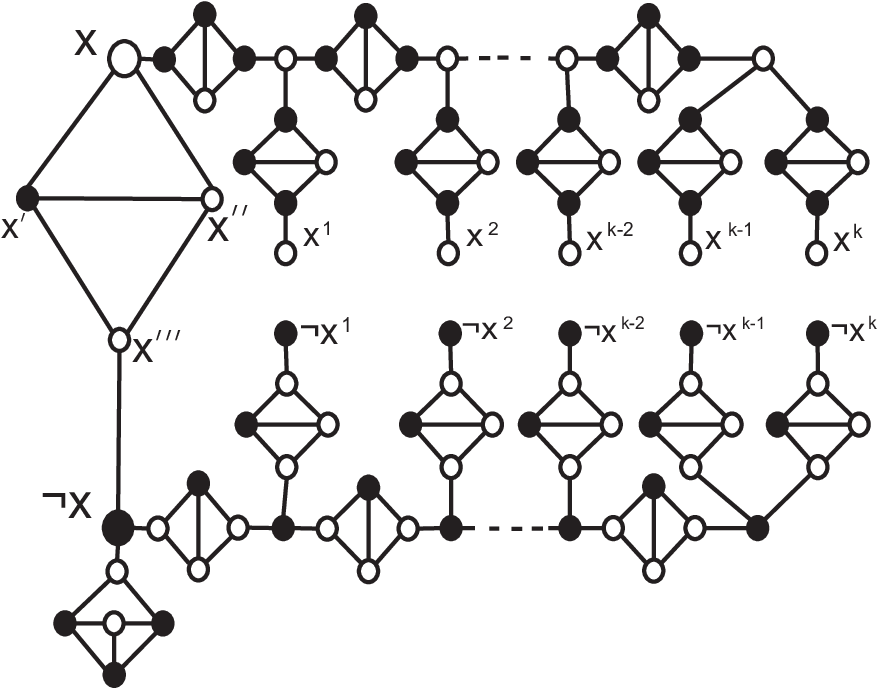}
				\caption{The graph $ H_x  $. In every lucky labeling $\ell: V(H_x) \rightarrow \mathbb{N}_2$, the set of black vertices has the same label and also the set of white vertices has the same label and these two labels are different (with respect to   symmetry).} \label{graphA}
			\end{center}
		\end{figure}
		
		Since for every regular graph $G$ we have that
		$\sigma(G)=2$ if and only if $\eta(G)=2$,  it suffices to prove the theorem for the lucky labeling.
		
		We reduce {\em NAE 3-SAT } to
		our problem in polynomial time.
		Consider an instance $ \Psi $ with the set of variables
		$X $ and   the set of clauses $C $. We transform this into a 3-regular graph $G_{\Psi}$
		such that $ \eta(G_{\Psi})=2$ if and only if $\Psi $ has an NAE truth assignment.
		We use three auxiliary graphs $ H_x$, $I_{c_j} $ and $T$. The gadgets   $ H_x$ and $T$ are shown in
		Figure \ref{graphA} and Figure \ref{graphB}. Also, the gadget $I_{c_j} $ is a triangle $ c_j^1 c_j^2 c_j^3 $. The graph $G_{\Psi}$ has
		a copy of $H_x $ for each variable $x \in X$ and a copy of
		$I_{c_j}$ for each clause $c_j \in C$.
		For each clause $c_j=y \vee z \vee w$, where  $y,w,z \in X \cup \neg X $, add the edges $c_j^1 y^j $, $c_j^2 z^j $ and $c_j^3 w^j $. Finally, for every vertex $v $ with $d(v)< 3$, put a copy of $T$ and add edge between $v$ and $t$. Repeat this procedure one more time to obtain a $3$-regular graph $G_{\Psi}$.

		\begin{figure}[ht]
			\begin{center}
				\includegraphics[scale=.6]{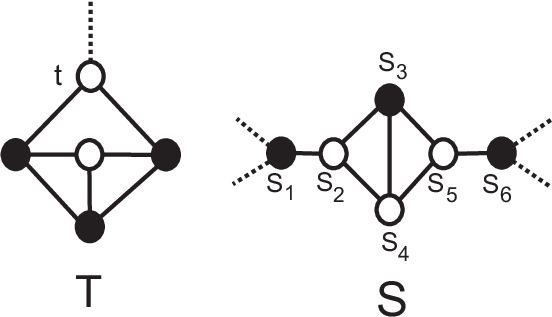}
				\caption{Two auxiliary graphs $ T$ and $S $. For every lucky labeling of $T$ with the labels $ \mathbb{N}_2$, the set of black vertices has a same label and also the set of white vertices has a same label and these two labels are different (with respect to   symmetry). This fact is also true for $S$.} \label{graphB}
			\end{center}
		\end{figure}
		
		We next discuss basic properties of the graph $G_{\Psi}$. Assume that $\eta(G_{\Psi})= 2$ and $\ell :V(G_{\Psi})\rightarrow \mathbb{N}_2$ is a lucky labeling.
		
		\begin{example}
			For every $x \in X$, in the subgraph $H_x$ we have, $\ell(x ) \neq \ell(\neg x )$.
		\end{example}
		{\bf Proof.}
		For a contradiction, suppose that  $\ell(x ) = \ell(\neg x )$. Consequently, $\ell \mid_{ \{ x' ,x'' ,x'''  \} }$ is a lucky labeling for the odd cycle $x' x'' x''' $. But the lucky number of every odd cycle is 3, a contradiction. $\spadesuit$
		
		If $u$ and $v$ are two vertices in the graph $G_{\Psi}$ such that $N[u]=N[v]$, then $\ell(u) \neq \ell(v)$. As a consequence, we have the following fact.
		
		\begin{example}
			In every copy of $S$ and $T$  (see Figure \ref{graphB}), the set of black vertices has the same label and also the set of white vertices has the same label and these two labels are different (with respect to   symmetry).
		\end{example}
		{\bf Proof.}
		In the subgraphs $S$, first assume that $\ell(s_1)=2$. Since $N[s_3]=N[s_4]$, we have that $\ell(s_3) \neq \ell(s_4)$. Without loss of generality suppose that $\ell(s_3)=2$ and $\ell(s_4)=1$. If $\ell(s_5)=2$, then $\ell \mid_{ \{s_2, s_3, s_4 \} }$ is a lucky labeling for the odd cycle $s_2 s_3 s_4$, but it is a contradiction, hence $\ell(s_5)=1$. If $\ell(s_2)=2$, then $ \sum_{w \sim s_2}\ell(w)= \sum_{w \sim s_4}\ell(w) $, so $\ell(s_2)=1$.
		Also if $\ell(s_6)=1$,
		then $ \sum_{w \sim s_4}\ell(w)= \sum_{w \sim s_5}\ell(w) $, so $\ell(s_6)=2$. By a similar argument if $\ell(s_1)=1$, then $\ell(s_6)=1$. Therefore, the set of black vertices has the same label and also the set of white vertices has the same label and these two labels are different.
		For the subgraph $T$ we have a similar argument.
		$\spadesuit$

		By Fact 1, Fact 2 and since $S$ is a subgraph of $H_x$, without loss of generality, for every $x \in X$, in the subgraph $H_x$, the set of black vertices has the same label and also the set of white vertices has the same label and these two labels are different. In other words, $\ell (x )=\ell(x^1)=\cdots = \ell(x^k) \neq \ell (\neg x)=\ell(\neg x^1)=\cdots = \ell(\neg x^k)$. For an arbitrary clause $c_j=y \vee z \vee w$,  where  $y,w,z \in X \cup \neg X $, assume that $  \ell( y^j)   =\ell(z^j)=\ell(   w^j )  $, consequently $\ell \mid_{ \{ c_j^1  ,  c_j^2 ,c_j^3   \} }$ is a lucky labeling for the odd cycle $ c_j^1   c_j^2 c_j^3   $, but the lucky number of an odd cycle is $ 3 $. This is a contradiction. Hence, we have the following fact:

		\begin{example}
			For every clause $c_j=y \vee z \vee w$, where $y ,z,w \in  X \cup \neg X $, we have $  \{ \ell(  y^j)   ,\ell(z^j),\ell(   w^j ) \} =\mathbb{N}_2$.
		\end{example}

		First, assume that $\eta(G_{\Psi})= 2$ and let $\ell :V(G_{\Psi})\rightarrow \mathbb{N}_2$ be a lucky labeling.
		We present an NAE
		satisfying assignment  $\Gamma : X\rightarrow \{\text{true}, \text{false}\}$ for $\Psi$. Now put $\Gamma(x_i)=true$ if and only if $\ell(x_i)=1$. By Fact 1, for every $x_i$ we have $\ell(x_i)\neq\ell(\neg x_i)$ so it is impossible that both $\ell(x_i)$ and $\ell(\neg x_i)$ are 1.
		For every $c_j=y \vee z \vee w$
		by Fact 3, $ \mid \{ \ell(  y^j)   ,\ell(z^j),\ell(   w^j ) \} \mid =2$; so at least one of the
		literals $y,z,w$ is true and at least one of the
		literals is false.
		On the other hand, suppose that $\Psi$ is satisfiable with the satisfying assignment $\Gamma : X  \rightarrow \{\text{true}, \text{false}\}$. We present the lucky labeling $\ell$ for $G_{\Psi}$ from the set $\{1,2\}$. By attention to  Figure \ref{graphA} and Figure \ref{graphB},  in each copy of $H_{x}$, if $\Gamma(x)= \text{true}$, label the set of black vertices by $1$ and label the set of white vertices by $2$. Also if $\Gamma(x)=\text{false}$, label the set of black vertices by $2$ and label the set of white vertices by $1$. For every copy of $T$, label the vertex $t$ with the label different from $\ell(t')$, where $t'$ is the unique neighbor of $t$ which is not in $T$. This determines the labels of all other vertices of $T$. By Fact $3$, for every $c_j\in C$,
		we can  determine the labels of the vertices $ c_j^1 ,c_j^2, c_j^3  $. Finally, by straightforward counting one can see that $\ell$ is a lucky
		labeling. This completes the proof.
}\end{ali}

\begin{alii}{
		Here, we prove that for every $k\geq 3$,  for a given graph
		$G$ it is
		{\bf NP}-complete to decide whether $\sigma(G)= k$. It was shown in  \cite{MR573644} that the following problem is NP-complete. We reduce it to our problem in polynomial time.
		\\ \\
		{\em The 3-colorability of 4-regular graphs.}\\
		\textsc{Instance}: A 4-regular graph $G$.\\
		\textsc{Question}: Is $\chi(G)\leq 3$?\\

		For a given $4$-regular graph $G$, we construct a regular graph $G^*$ such that $\chi(G^*)\leq k$  if and only if
		$G$ is 3-colorable (step 1). Next, we construct the graph $G^{**}$ such that $\sigma(G^{**})\leq k$ if and only if $\chi(G^*)\leq k$ (step 2).
		
		Step 1. If $k=3$ we can assume that $G^* \cong G$, otherwise do the following. Consider a copy of the graph $G$ and a copy of the complete graph $K_{k-3}$. Join each vertex of $G$ to each vertex of $K_{k-3}$ and call the resultant graph $G'$. One can see that $\chi(G')=\chi(G)+k-3$  and $\Delta(G')=k-4+|V(G)|$. Now, consider a copy of the graph $G'$ and for each vertex $v\in G'$, put $\Delta(G')-d_{G'}(v)$ new isolated vertices and join them to the vertex $v$. Also, put $\Delta(G')$ copies of the $K_2$'s.
		Call the resultant graph $G''$. It is easy to check that $\delta(G'')=1$, $\Delta(G'')=\Delta(G')$ and $\chi(G'')=\chi(G')$.
		Let $S=\{v \mid d_{G''}(v)=1\}$.
		Now, consider two copies of $G''$ and put $\Delta(G'')-1$ distinct perfect matching between the set of vertices $S$ in the first copy of $G''$ and the set of vertices $S$ in the second copy of $G''$. Call the resultant graph $G^*$. On can see that $G^*$ is a regular with $\chi(G^*)=\chi(G)+k-3$.

		Step 2. Suppose that $G^*$ is a regular graph with $n$ vertices. For every $\alpha$, $ 1 \leq \alpha \leq n $ consider a copy of a complete graph $K^{(\alpha)}_{k^2-k+1}$, with the vertices $\{  y_{\beta \gamma }^{\alpha} :   \beta, \gamma \in \mathbb{N}_{k-1}  \}  \cup \{ y_{k \gamma }^{\alpha}:    \gamma \in \mathbb{N}_{k }   \}$. Next consider $k-1$ isolated vertices $v_1, \ldots , v_{k-1}$ and join every $y_{\beta \gamma }^{\alpha}$ to $v_{\beta}$,
		also consider a copy of $G^*$ with the vertices $x_1, \ldots , x_{n}$ and join every $y_{k \gamma }^{\alpha}$ to $x_{\alpha}$. Finally put $n$ isolated vertices $z_1, \ldots , z_{n}$ and join every $x_{\alpha}$ to $z_{\alpha}$. Name the constructed graph $G^{**}$.
		
		First, note that $\sigma(G^{**})\geq k$, since in every sigma partitioning of $G^{**}$ the labels of $y^1_{k1},y^1_{k2},\ldots, y^1_{kk}$ are different. Let $c :V(G^{**}) \rightarrow \{a_1, \ldots, a_k \}$ be a sigma partitioning. For every $\alpha$ and $\beta$, the labels of $y^{\alpha}_{\beta 1}, \cdots, y^{\alpha}_{\beta (k-1)}$ are different.
		We claim that the labels of  $v_1, \ldots , v_{k-1} $ are different (Property  A). To the contrary suppose that $c(v_{\beta})=c(v_{\beta '})$, then the labels of the vertices in $\{ y_{\beta \gamma }^{1}:   \gamma\in \mathbb{N}_{k-1} \} \cup \{  y_{\beta \gamma' }^{1}:   \gamma' \in \mathbb{N}_{k-1} \}$ must be different, so we need at least $2(k-1)$ labels. Similarly, for every $i$ and $\gamma$, $  i  \in \mathbb{N}_{n} $, $   \gamma \in \mathbb{N}_{k-1}$, we can see that $ c(v_{\beta}) \neq c(x_i)$ (Property  B). First, suppose that $\sigma(G^{**})= k$, by Properties  A and B,
		we have $c(x_1)=\cdots =c(x_n)$. Consider the set of vertices $\{x_i : i \in \mathbb{N}_{n}\}$, since $G^*$ is regular and for every $\alpha$ we have $|\{ c(y^{\alpha}_{ k \gamma}):    \gamma \in \mathbb{N}_k   \}|=k$, therefore, the function
		$c':\{x_i : i \in \mathbb{N}_{n}\}\rightarrow \mathbb{N} $, where $c'(x_{\alpha})= c(z_{\alpha})$, is a proper vertex coloring of $G^*$,
		so $\chi(G^*)\leq k$. On the other hand, let $\chi(G^*)\leq k$ and $c':\{x_i : i \in \mathbb{N}_{n}\} \rightarrow \mathbb{N}_k$ be a proper vertex coloring $G^*$. Suppose that $s$ is a sufficiently large number, we give a sigma partitioning $c$ for $G^{**}$ with the labels $\{s^i : i \in \mathbb{N}_{k} \}$, that is:
		
		$c(x_{\alpha})=s^k$, $ $ $ $ $ $  $ $ $c(v_{\beta})=s^{\beta}$, $ $ $ $ $ $ $ $
		$c(y^{\alpha}_{k \gamma})=s^{\gamma}$, $ $ $ $  $ $  $ $ $c(z_{\alpha})=s^{c'(x_{\alpha})}$,
		\\ \\
		and for every $\alpha$, $\beta$ with $\beta \neq k$,  label the vertices $y^{\alpha}_{ \beta \gamma}$ such that $ \{c(y^{\alpha}_{ \beta \gamma}):  \gamma \in \mathbb{N}_k \}=\{ s^i:   i \in \mathbb{N}_k \}  \setminus \{ s^{\beta}\}$.
		It is not hard to check that $c$ is a sigma partitioning of $G^{**}$. This completes the proof.
}\end{alii}

\begin{aliii}{
		Here, we prove that for a given planar $3$-regular graph $G$ with  $\sigma(G) = 2$ and a real number $0 < r  <  1$,
		determining whether $G$ has a sigma partitioning with $\sigma(G)$ parts, such that there is a part with at most $r|V(G)|$ vertices, is $ \mathbf{NP} $-complete.
		\\
		First consider the following problem.
		\\ \\
		{\em Planar NAE 3-SAT.}\\
		\textsc{Instance}: Set $X$ of variables, collection $C$ of clauses over $X$ such that each
		clause $c \in C$ has $| c |= 3$ and the following graph obtained from 3-SAT is planar. The graph has one vertex for each variable, one vertex for each clause; all variable vertices are connected
		in a simple cycle and each clause vertex is connected by an edge to variable
		vertices corresponding to the literals present in the clause.\\
		\textsc{Question}: Is there an NAE truth assignment for $X$?\\

		Moret \cite{p} proved  that  {\em Planar NAE 3-SAT} is in $ \mathbf{P} $ by a reduction to a
		known problem in $ \mathbf{P} $, namely Planar MaxCut (for more   information see \cite{MR3544062, MR3864719}). Moret's reduction used
		only local replacement. Given an instance of {\em Planar NAE 3-SAT} with $k$ clauses, they transformed it in polynomial time into an instance $Q$ of MaxCut with $9k$ vertices as follows.
		For each variable $x$ forming a total of $n_x$ literals in the $k$ clauses, they put
		a cycle of $2n_x$ vertices. Alternating vertices represent
		complemented and uncomplemented literals. For each clause, they put a copy of the complete graph $K_3$, where each vertex of the $K_3$ is connected by an edge to the (complement of the) corresponding literal vertex. They proved that $11k$ is the maximum attainable
		cut sum if and only if {\em Planar NAE 3-SAT} is satisfied. Note that the cycle that contains all variables does not have any rule in the construction of $Q$. So, if we replace an instance of  {\em Planar NAE 3-SAT} with an instance of the following problem; the proof remains correct. Therefore, the following problem is in  $ \mathbf{P} $.
		\\ \\
		{\em Planar NAE 3-SAT Type 2.}\\
		\textsc{Instance}: Let $X$ be the set of variables and let $C$ be
		the set of clauses such that each
		clause $c \in C$ has $|c  | = 3$ and the
		bipartite graph obtained by linking a variable and a clause if and
		only if the variable appears in the clause, is planar.\\
		\textsc{Question}: Is there an NAE truth assignment for $X$?\\

		On the other hand, Moore and Robson \cite{MR1863810} proved
		that the following 1-In-3 SAT problem is $ \mathbf{NP}$-complete.
		\\ \\
		{\em  Cubic Planar 1-In-3 3-SAT.}\\
		\textsc{Instance}: Set $X$ of variables, collection $C$ of clauses over $X$ such that each
		clause $c \in C$ has $| c  | = 3$ and every variable appears in
		exactly three clauses, there is no negation in the formula, and the
		bipartite graph obtained by linking a variable and a clause if and
		only if the variable appears in the clause, is planar. Note that $|C| =|X|$.\\
		\textsc{Question}: Is there a truth assignment for $X$ such that each clause in $C$ has exactly
		one true literal?\\

		\begin{figure}[ht]
			\begin{center}
				\includegraphics[scale=.6]{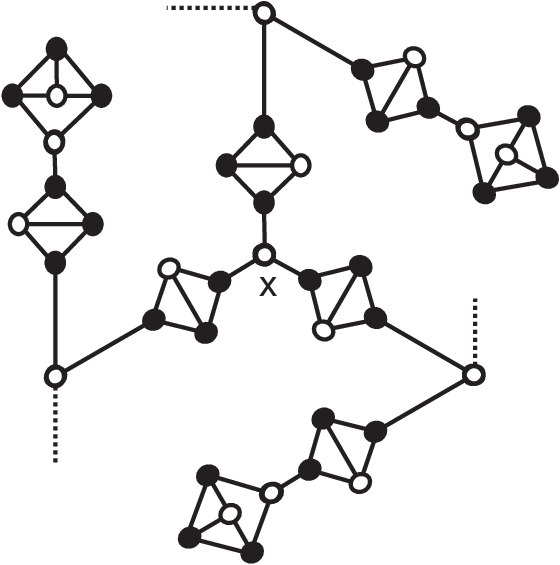}
				\caption{The auxiliary graph $A_x$. For every sigma partitioning of $A_x$ with the labels $ \{ \alpha,\beta\}$, the set of black vertices have the same label and also the set of white vertices have the same label and these two labels are different (with respect to the symmetry).
				} \label{graphc}
			\end{center}
		\end{figure}
		
		We reduce {\em Cubic Planar 1-In-3 3-SAT} to
		our problem in polynomial time.
		Assume that $H $ is a dummy 3-regular graph with
		$\sigma(H ) = 2$. Also, without loss of generality, suppose that $H $ does not have any sigma partitioning with two parts, such that there is a part with at most $\frac{62}{138}|V(H)|$ vertices.
		Now, consider an instance $ \Phi $ with the set of $3n$ variables
		$X $ and the set of $3n$ clauses $C  $.
		By the polynomial time algorithm which was presented in \cite{p} (they proved that {\em Planar NAE 3-SAT Type 2} is in $ \mathbf{P} $), we check whether $\Phi$ has an NAE truth assignment. If $\Phi$ does not have any NAE truth assignment, then it does not have a 1-In-3 truth assignment. In this case we return $H$.  Otherwise, we construct $G_{\Phi}$.
		We use two auxiliary graphs $ A_x$ and $ B_c $. The gadget $ A_x$ is shown in
		Figure \ref{graphc} and the gadget $B_c$ is a complete graph $K_3$. Let $H_{\Phi}$ be the
		graph obtained by linking a variable and a clause if and
		only if the variable appears in the clause. Replace each variable $x \in X$ by $A_x$ and replace each clause $c \in C$ by $B_c$. Call this planar $3$-regular graph $G_{\Phi}$. This graph has $138n$ vertices.
		
		First we show that $\sigma(G_{\Phi}) = 2$. By assumption $\Phi$ has an NAE truth assignment $\Gamma : X \rightarrow \{true, false\}$. We present a sigma partitioning $c$ for $G_{\Phi}$, with labels $\alpha$ and $\beta$.
		Put $c(x)=\beta$ if and only if $\Gamma(x)=true$.
		This determines the labels of the remaining vertices in $A_x$.
		Let $c\in C$ be an arbitrary clause and without loss of generality, assume that $ v_1v_2v_3$ is the complete graph corresponds to the clause $c=(x \vee y \vee z)$ in the graph. Also, assume that $v_1$ has a neighbor in $A_x$ with label $\alpha$,
		$v_2$ has a neighbor in $A_y$ with label $\beta$ and $v_3$ has a neighbor in $A_z$ with label $\beta$. Put $c(v_2)=\alpha$ and
		$c(v_1)=c(v_3)=\beta$.
		Since $\Gamma$ is an NAE truth assignment, one can  label the vertices of $B_c$ by the above method.

		We prove that for every sigma partitioning with $\sigma(G_{\Phi})$ labels, for every part $P$ we have $|P|/|V(G_{\Phi})|\geq \frac{62}{138}$ and the equality holds for a part,
		if and only if $\Phi$ has a truth assignment such that each clause in $C$ has exactly
		one true literal. Assume that  $\ell :V(G_{\Phi})\rightarrow \{\alpha , \beta\}$ is a sigma partitioning.
		Consider the following fact:

		\begin{example}\label{F4}
			Let $G$ be a regular graph, $\sigma(G)=2$ and $\ell :V(G)\rightarrow \{\alpha , \beta\}$  be a sigma partitioning. Also, let
			$\alpha' $ and $\beta'$ be two arbitrary number such that $\alpha' \neq\beta'$. Define $\ell :V(G)\rightarrow \{\alpha' , \beta'\}$ such that $\ell(v)=\alpha'$ if and only if $\ell(v)=\alpha$. The function $\ell$ is a sigma partitioning for the graph $G$.
		\end{example}
		{\bf Proof.}
		Since the graph $G$ is regular, the proof is clear.
		$\spadesuit$
		
		By Fact \ref{F4} and Fact 2,  in the subgraph $A_x$, the label of $x$ forces the labels of all other vertices; therefore,
		the labels of exactly $16$ vertices including $x$
		are equal to $\ell(x)$. Also, by the structure of $B_c$, the function $\Gamma$ is an NAE satisfying  assignment. So, in every $B_c$, at least one vertex has label $\beta$, at least one vertex has label $\alpha$, also at least one vertex of their neighbors in $ \cup_{x \in X} A_x$  has label $\beta$  and  at least one vertex of their neighbors in $ \cup_{x \in X} A_x$  has label $\alpha$.
		As $|C| =|X|=3n$ and $\Phi$ is satisfied, at
		least $\frac{1}{3}$ of the variables are true and at
		most $\frac{2}{3}$ of the variables are false. Thus  for $\alpha$, we have:
		
		\begin{align*}
		|  \ell^{-1}(\alpha) | &= |  \ell^{-1}(\alpha) \cap \bigcup_c B_c | + |  \ell^{-1}(\alpha) \cap \bigcup_{\ell(x)=\alpha} A_x|+ |  \ell^{-1}(\alpha) \cap \bigcup_{\ell(x)=\beta} A_x| \\
		&\geq (1 \times 3n)+(16 \times 2n )+ (27 \times n)\\
		&=62n
		\end{align*}

		If the  equality holds, then $\Phi$ has a 1-In-3
		satisfying assignment with the satisfying assignment $\Gamma : X \rightarrow \{true, false\}$, when $\Gamma(x)=true$ if and only is $\ell(x)=\beta$.
		Next, suppose that $\Phi$ has 1-In-3
		satisfying assignment with the satisfying assignment $\Gamma : X \rightarrow \{true, false\}$.
		We present the sigma partitioning $\sigma$ for $G_{\Phi}$ from the set $\{\alpha,\beta\}$ such that exactly $62n$ vertices have the label $\beta$. Put $\ell(x)=\beta$ if and only if $\Gamma(x)=true$, the labels of other vertices of $A_x$ are determined.
		Now, for each $c \in C$, exactly one of the vertices of $B_c$ has a neighbor in $\cup_x A_x$ with the label $\beta$. Call this vertex of $B_c$ by $v_c$. Label $v_c$ by $\alpha$ and label the two  other vertices of $B_c$ by different labels. One can see that $\sigma$ is a sigma partitioning.  This completes the proof.
}\end{aliii}

\section {Random graphs and upper bounds}

\begin{aliiii}{
		(i) Let $\mathcal{G} = \{G_i\}_{i \in \mathbb{N}}$ be a sequence  of triangle-free graphs. Then for this family we prove that if $\delta =  \Omega(\Delta)$ then $\sigma =  \mathcal{O} (1) $, for all graphs in $\mathcal{G}$, except a finite number of them.
		We will use the probabilistic method to prove the theorem. The following tool
		of the probabilistic method will be used several times.
		
		\begin{lam}
			\noindent {\bf [The Local Lemma \rm \cite{MR1885388}\bf]} Suppose ${A_{1},\ldots ,A_{n}}$ is a set of random events such that for each
			$i$, $\Pr(A_{i})\leq p$ and $A_{i}$ is mutually independent of the set of all but at most $d$ other events. If $ep(d+1)< 1$, then with positive probability, none of the events occur.
		\end{lam}
		
		We will also use the following well known inequalities (Stirling's approximation).
		
		$$\sqrt{2\pi n} (\frac{n}{e})^n \leq n! \leq \sqrt{e^2 n} (\frac{n}{e})^n.$$
		
		Let $G$ be a given graph. 
		We will use the Local Lemma in order to prove the Theorem.
		Color each vertex of $G$ by a random color from the set $\mathbb{N}_k$
		such that each color is chosen with probability $\frac{1}{k}$
		and the color of every vertex is independent from the colors of other vertices.
		The value of $k$ will be determined later. For an edge $e=uv$,
		let $B_e$ be the event {\it ``for each $i$, $i \in \mathbb{N}_k$, the vertex $u$ and the vertex $v$ have
			the same number of neighbors in part $i$"}.
		Note that if $B_e$ occurs, then $d(u)=d(v)$.
		For every edge $e$, the event $B_e$ is dependent on less than $2 \Delta ^3$ events.
		Since the maximum of $\Pr(B_e)$ is when $d(v)=d(u)=\delta(G)$ and $G$ is triangle-free, we have:

		\begin{align}
		\Pr(B_e)  & \leq \displaystyle\sum_ {\sum_ {i} t_i=\delta}\Big(\dfrac{\displaystyle{{\delta}\choose {t_1, \ldots, t_k}}}{k^ \delta} \Big)^2\\
		& \leq \displaystyle \max_{\sum_{i} t_i=\delta} \dfrac{\displaystyle{{\delta}\choose{t_1, \ldots, t_k}} }{k^ {\delta}} \\
		& \leq  \dfrac{\displaystyle{{\delta}   \choose   {\frac{\delta}{k} , \ldots, \frac{\delta}{k}}   }}{ k^ {\delta}}\\
		&= \dfrac{\delta!}{\Big((\frac{\delta}{k})!\Big)^k}\times k^{-\delta}
		\end{align}       
		By Stirling's approximation and (4), we have       
		\begin{align}         
		\Pr(B_e)   &\leq  \frac{e\sqrt{\delta}(\delta/e)^ \delta}{\Big(\sqrt{\frac{2\pi \delta}{k}}(\frac{\delta}{ek})^ {\frac{\delta}{k}} \Big)^k}\times k^{-\delta} \\
		&  = \frac{e\sqrt{\delta}}{\Big (\sqrt{\frac{2\pi \delta}{k}}\Big)^k}.
		\end{align}
		It was shown in \cite{survey} that  for every graph $G$, we have $\sigma(G)= \mathcal{O} (\Delta^2)$.
		So for a given graph $G$ if $\Delta(G)$ is a constant number then $\sigma(G)=\mathcal{O} (1)$.
		Thus, we can assume that for a given graph $G$,  $\Delta$ is a sufficiently large number.
		Since $\delta =  \Omega(\Delta)$, there are constant numbers $n_0$ and $c$ such that for each $i\geq n_0$, for the graph $G_i$,
		we have $\delta \geq c \Delta $. Put $k=8$. For each $i\geq n_0$, for the graph $G_i$ with a sufficiently large $\Delta$, we have $ ep(d+1)\leq e \frac{e\sqrt{\delta}}{\Big (\sqrt{\frac{2\pi \delta}{8}}\Big)^8} \Delta^3<1$. Thus by the Local Lemma there is a sigma partitioning with $\mathcal{O} (1)$ colors   for all graphs in $\mathcal{G}$, except a finite number of them.
		This completes the proof of the part (i).
		\\ \\
		(ii) Let $\mathcal{G} = \{G_i\}_{i \in \mathbb{N}}$ be a sequence  of triangle-free graphs. Then for this family we prove that if $\delta \geq  3.501 (\ln \Delta)$, then $\sigma = \mathcal{O}(\delta) $.
		Put $k=\delta$. We use the Local Lemma to prove the theorem. From the previous calculation we have:
		
		\begin{align}
		\Pr(B_e) & \leq  \dfrac{\displaystyle{{\delta}   \choose   {\frac{\delta}{k} , \ldots, \frac{\delta}{k}}   }}{ k^ {\delta}}
		\end{align}
		Since $k=\delta$, by (7), we have
		\begin{align}
		\Pr(B_e) & \leq       \dfrac{\displaystyle{{\delta}   \choose   {  1 , \ldots, 1                            }   }}{ \delta^ {\delta}}  \\
		& = \frac{\delta !}{\delta^{\delta}}
		\end{align}         
		By Stirling's approximation and (9), we have       
		\begin{align}         
		\Pr(B_e)         &\leq \frac{e\sqrt{\delta}}{e^{\delta }}
		\end{align}
		It was shown in \cite{survey} that  for every graph $G$, we have $\sigma(G)= \mathcal{O} (\Delta^2)$.
		So for a given graph $G$ if $\Delta(G)$ is a constant number then $\sigma(G)=\mathcal{O} (1)$.
		Thus, we can assume that for a given graph $G$,  $\Delta$ is a sufficiently large number.
		Since $\delta \geq  3.501 (\ln \Delta)$ and $\Delta $ is a sufficiently large number, we have $ep(d+1)\leq e \frac{e\sqrt{\delta}}{\Delta^{3.501}} \Delta^3  < 1$. Thus by the Local Lemma there is a sigma partitioning. This completes the proof.
		
}\end{aliiii}

\begin{aliiiii}{
		Here, we prove that for every constant $p$, $0<p<1$, $\sigma(G(n,p))\leq 5$.
		The proof is similar to the proof of
		Theorem 4 in \cite{MR3072733}.
		Suppose that $\{\mathcal{P}_1, \ldots , \mathcal{P}_5 \}$
		is a partition for $n$ vertices such that
		the size of each part is $\lfloor n/5 \rfloor $ or $ \lceil n/5 \rceil$.
		Next consider the graph $G(n, p)$ on these vertices.
		For every vertex $v$ denote the number of
		neighbors of the vertex $v$ in $\mathcal{P}_i$ by $\mathcal{A}_v^i$.
		This partition is a sigma partitioning if for every
		two adjacent vertices $v$ and $u$, there is an
		index $i$, such that $\mathcal{A}_v^i \neq \mathcal{A}_u^i$. First, we calculate $\Pr(\mathcal{A}_v^i = \mathcal{A}_u^i) $. There is a constant $C$ such that 
		\begin{equation}
    	\displaystyle \Pr(\mathcal{A}_v^i = \mathcal{A}_u^i)
       \leq C\Big( \displaystyle\sum_{t=0}^{  \frac{n}{5}  } \big( { {\frac{n}{5}} \choose   {t}  }p^t (1-p)^{  \frac{n}{5} -t }\big)^2 \Big) 
        \end{equation}
        By (11), there is a constant $C'$ such that
        \begin{equation}
        \displaystyle \Pr(\mathcal{A}_v^i = \mathcal{A}_u^i)
        \leq C'  \Big(\max_{0\leq t \leq  n/5  }{ {\frac{n}{5}} \choose   {t}  }p^t (1-p)^{  \frac{n}{5} -t }\Big) 
        \end{equation}
        On the other hand, by Stirling's approximation we have
         \begin{equation}
        \Theta \Big(\max_{0\leq t \leq  n/5  }{ {\frac{n}{5}} \choose   {t}  }p^t (1-p)^{  \frac{n}{5} -t }\Big) =\Theta \Big( { {\frac{n}{5}} \choose   {\frac{pn}{5}}  } p^{\frac{p  n}{5}   } (1-p)^{\frac{(1-p)  n}{5}}\Big) = \Theta(n^{-\frac{1}{2}})
        \end{equation}
      Thus, by (12) and (13), we have $\Pr(\forall i \ \mathcal{A}_v^i = \mathcal{A}_u^i) =  \Theta(n^{-\frac{5}{2}})$, therefore
		$ \Pr(\exists vu \ \forall i \ \mathcal{A}_v^i = \mathcal{A}_u^i) = \Theta(n^2) \Theta(n^{-\frac{5}{2}})=o(1)$.
		This completes the proof.
}\end{aliiiii}

\section{Concluding remarks}
\
\textbullet $ $ A hypergraph $\mathcal{H}$ is a pair $(X,Y)$,
where $X$ is the set of vertices and $Y$ is a set of non-empty subsets of $X$, called edges. A $k$-coloring of $\mathcal{H}$ is a coloring $\ell:X \rightarrow \mathbb{N}_{k}$ such
that, for every edge $e$ with $| e| > 1$, there exist $v,u\in e$ such that $\ell(u) \neq \ell(v)$ (for more information about hypergraphs, see \cite{graham1995handbook}).
We say that two vertices $v$ and $u$ are adjacent if there is an edge like $e$ such that $v,u\in e$.
A sigma partitioning for a hypergraph $\mathcal{H}$
is a partition for the vertices
such that every edge $e$ has two
adjacent vertices $v$ and $u$, such that $v$ and $u$ have different numbers of neighbors
in some parts. Investigating the property of sigma partitioning in hypergraphs can be interesting as a future work.
\\
\textbullet $ $ In Theorem \ref{T1}, we used {\em NAE 3-SAT }. The planar version of {\em NAE 3-SAT }
is in $ \mathbf{P} $ \cite{p}, so the computational complexity of $\sigma $ for planar 3-regular graphs remains unsolved.

\bibliographystyle{plain}
\bibliography{luckyref}

\end{document}